\documentclass[copyright,creativecommons noncommercial]{eptcs}

\usepackage{iftex}
\ifpdf
  \usepackage{underscore}         
  \usepackage[T1]{fontenc}        
\else
  \usepackage{breakurl}           
\fi

\usepackage{amsmath,amsfonts}
\usepackage{graphicx}
\usepackage{hyperref}
\usepackage{caption}
\usepackage{subcaption}

\newcommand{\R}{\mathbb{R}}

\title{HBS Tilings Extended:\\ 
State of the Art and Novel Observations}
\author{Carole Porrier
\institute{Université du Québec à Montréal, Canada}
\institute{Université Sorbonne Paris Nord, France}}
\def\titlerunning{HBS Tilings Extended}
\def\authorrunning{C. Porrier}

\begin{document}

\maketitle

\begin{figure}[h!]
    \centering
    \begin{subfigure}[b]{0.35\linewidth}
        \centering
        \includegraphics[scale=0.25]{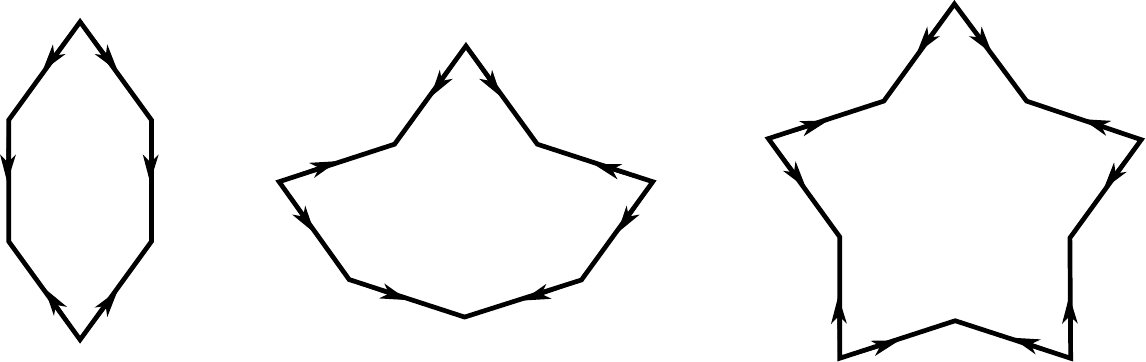}
        \caption{HBS tiles (Hexagon-Boat-Star).}\label{fig:hbs-tiles}
    \end{subfigure}
    \hfill
    \begin{subfigure}[b]{0.6\linewidth}
        \centering
        \includegraphics[scale=0.072]{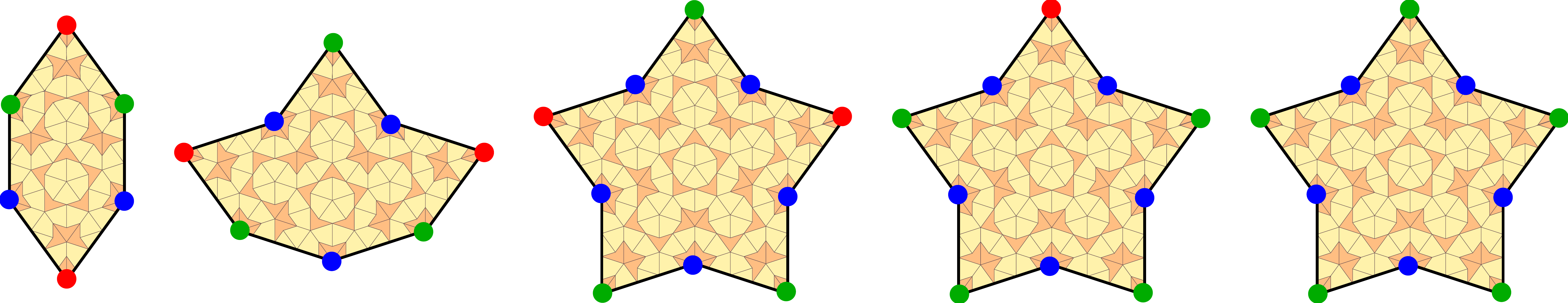}
    \caption{Star tileset with small kites and darts decorations.}\label{star-tileset-kd}
    \end{subfigure}
    \caption{Usual HBS tiles and their new, enriched version.}
\end{figure}

Penrose tilings are the most famous aperiodic tilings since Gardner described them \cite{gardner1977}, and they have been studied extensively \cite{grunbaum1987,Senechal1995,BaakeGrimm2013}.
It is thus surprising that one can still discover something new about them, even if it is not completely new.
Indeed, patterns composed with hexagons $H$, boats $B$ and stars $S$ (HBS tiles, Fig. \ref{fig:hbs-tiles}) were soon exhibited and aroused interest among physicists.
As serendipity works, we rediscovered HBS shapes but with different decorations (Figure \ref{star-tileset-kd}) and forcing rules, while working on a combinatorial optimization problem on graphs defined by kites-and-darts Penrose tilings (type P2), as described in \cite{porrier2020}.
Additionally, we distinguished three types of P2-stars depending on their surrounding, and as we labeled the HBS-vertices according to these types it appeared that most vertices of a given shape always had the same label, with an exception for the star.
We thus called \textit{Star tileset} the five shapes with labeled vertices, before knowing about the existence of HBS tiles -- even though they first appeared in Henley's 1986 paper \cite{Henley1986}.

As we found it difficult to know whether some of our findings had already been published and if so, including which information, it seemed a state of the art would be useful to us and might be for others too.
But most of all, we have new findings which we rely on for an article to come \cite{porrier2023}.


\section{Penrose tilings and mutual local derivability}

A \textit{tiling} of $\R^2$ is a countable family of non-empty closed sets $(T_i)_{i\in I}$ called \textit{tiles} such that $\bigcup_{i\in I}T_i=\R^2$ and $\mathring{T_i}\cap\mathring{T_j}=\emptyset$ for all $i\neq j$ in $I$.
The \textit{prototiles} of a tiling are the equivalence classes of its tiles up to congruence (in the present context).
A set of prototiles is called a \textit{tileset}, and oftentimes many tilings can be composed with copies of the same prototiles.
The three types of Penrose tilings, denoted P1, P2 and P3 in \cite{grunbaum1987}, were described by Roger Penrose himself \cite{penrose1978}. 
\begin{figure}[h!]
    \centering
    \begin{subfigure}[b]{0.5\columnwidth}
        \centering
        \includegraphics[width=\linewidth]{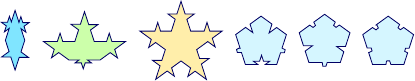}
        \caption{Tiny diamond, boat, star and pentagons (P1)}\label{P1}
    \end{subfigure}
    \hfill
    \begin{subfigure}[b]{0.2\columnwidth}
        \centering
        \includegraphics[scale=0.4]{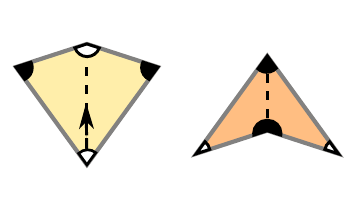}
        \caption{Kite and dart (P2)}\label{P2}
    \end{subfigure}
    \hfill
    \begin{subfigure}[b]{0.22\columnwidth}
        \centering
        \includegraphics[scale=0.4]{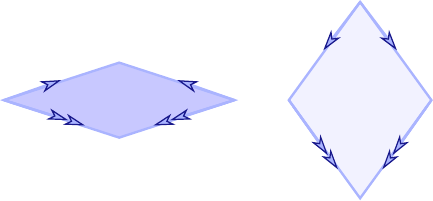}
        \caption{Rhombuses (P3)}\label{P3}
    \end{subfigure}
    \caption{Penrose tilesets.}
    \label{fig:Penrose-tilesets}
\end{figure}
The corresponding tilesets are shown in Figure \ref{fig:Penrose-tilesets}. 
Uncountably many tilings can be composed with each, but none of them is periodic: they have no translations among their symmetries.
The tiles must be arranged according to specific assembly rules, for instance using notches and bumps to assemble the tiles like puzzle pieces (P1 tiles, Fig. \ref{P1}) or markings on the tiles which can be of different kinds.
For P2, we use two colors which must match on the corners of the tiles (Fig. \ref{P2}), forming a full black or blank disk at each vertex of a tiling. 
Such a marking would not be sufficient for P3 tiles so we use arrows on the edges (Fig. \ref{P3}), which must superimpose exactly (same number of arrows, in the same direction).
This marking is also the most convenient in this article, in relation with HBS tiles.
For the same further purposes, an additional arrow is added on the diagonal of the kite, which corresponds to the single arrow on the edges of rhombuses when a P2 tiling is derived from a P3 tiling.
The sides of rhombuses have length 1, which is also the length of the kite's diagonal. 
Hence for kites and darts the longer side has length 1 and the shorter $\varphi^{-1}$ where $\varphi = (1+\sqrt{5})/2$ is the golden ratio. 
There are only seven vertex configurations, that is seven ways in which kites and darts can be arranged around a vertex of the tiling (Figure \ref{vertex-config}).
%
\begin{figure*}[h]
    \centering
    \begin{subfigure}[b]{0.11\textwidth}
        \centering
        \includegraphics[scale=0.3]{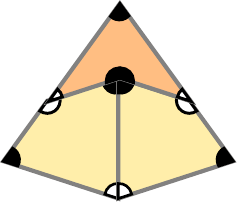}
        \caption*{Ace}
        \label{ace}
    \end{subfigure}
    \hfill
    \begin{subfigure}[b]{0.11\textwidth}
        \centering
        \includegraphics[scale=0.3]{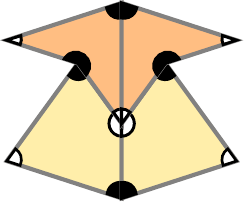}
        \caption*{Deuce}
        \label{deuce}
    \end{subfigure}
    \hfill
    \begin{subfigure}[b]{0.11\textwidth}
        \centering
        \includegraphics[scale=0.3]{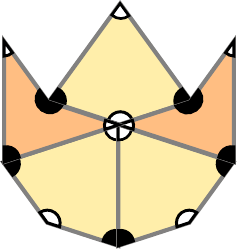}
        \caption*{Jack}
        \label{jack}
    \end{subfigure}
    \hfill
    \begin{subfigure}[b]{0.11\textwidth}
        \centering
        \includegraphics[scale=0.3]{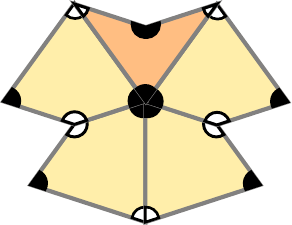}
        \caption*{Queen}
        \label{queen}
    \end{subfigure}
    \hfill
    \begin{subfigure}[b]{0.11\textwidth}
        \centering
        \includegraphics[scale=0.3]{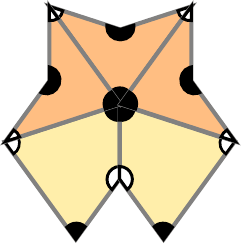}
        \caption*{King}
        \label{king}
    \end{subfigure}
    \hfill
    \begin{subfigure}[b]{0.11\textwidth}
        \centering
        \includegraphics[scale=0.3]{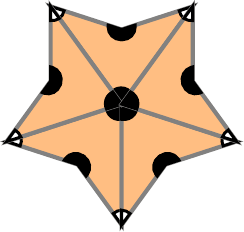}
        \caption*{Star}
        \label{star}
    \end{subfigure}
    \hfill
    \begin{subfigure}[b]{0.11\textwidth}
        \centering
        \includegraphics[scale=0.3]{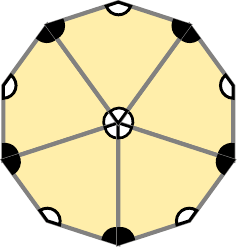}
        \caption*{Sun}
        \label{sun}
    \end{subfigure}
    \caption{The seven vertex configurations in a Penrose tiling by kites and darts.}
    \label{vertex-config}
\end{figure*}

Any Penrose tiling (of any type) can be composed or decomposed into another Penrose tiling, of the same or another type, using a local mapping: those tilings are \textit{mutually locally derivable} (\textit{MLD}).
For instance, Penrose rhombs can be decomposed into kites and darts, which can in turn be decomposed into smaller rhombs.
An \textit{inflation} (in the case of Penrose tilings) is a decomposition into tiles of the same type as the original ones, followed with a $\varphi:1$ scaling, so that the new tiles have the same shapes and sizes as the original ones.


\section{HBS tilings summarized}\label{sec:hbs-stateoftheart}

\begin{figure}[h!]
    \centering
    \begin{subfigure}[b]{0.48\textwidth}
    \centering
    \includegraphics[scale=.25]{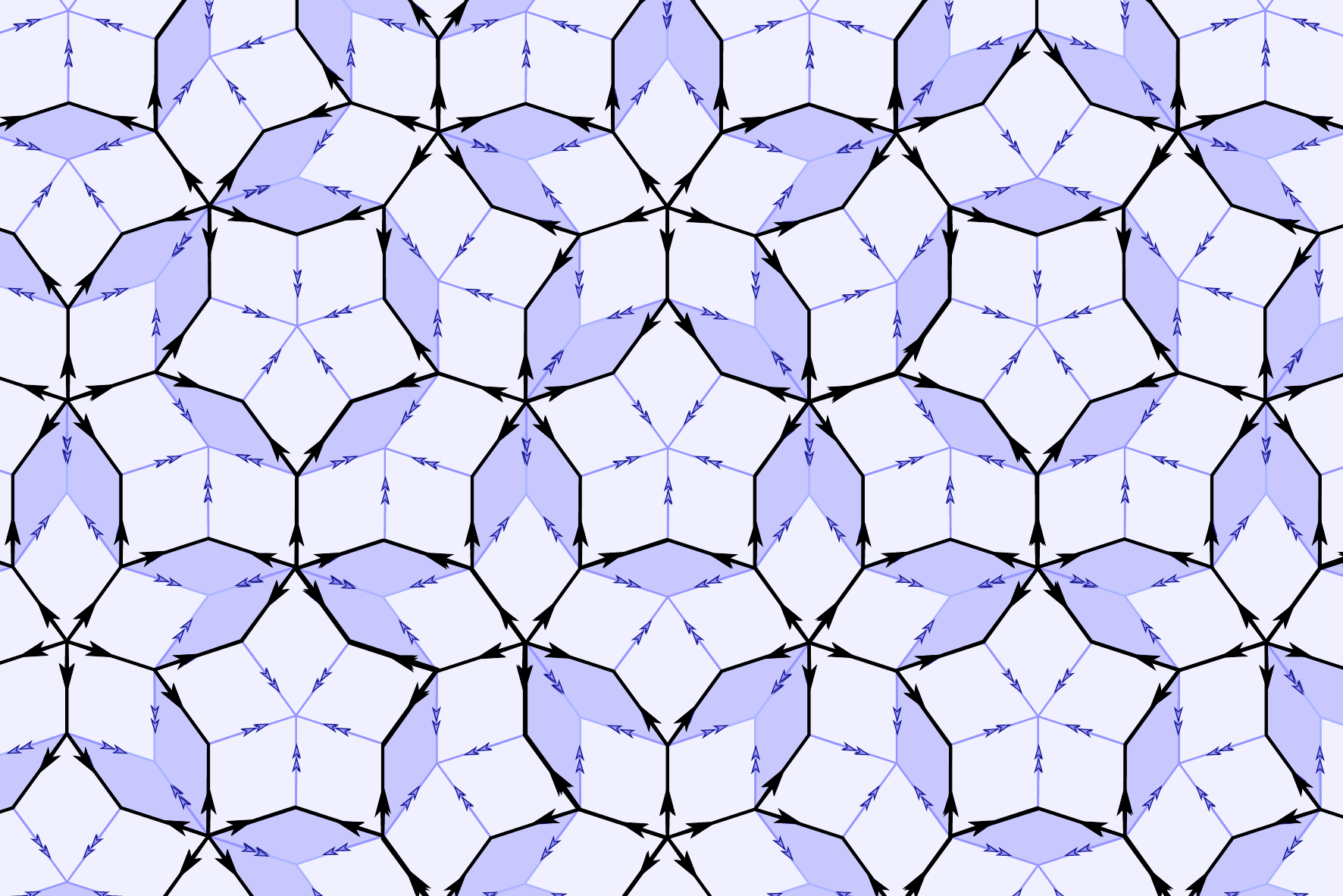}
    \caption{HBS and P3 Penrose tilings superimposed.}\label{fig:hbs-P3}
    \end{subfigure}
    \hfill
    \begin{subfigure}[b]{0.48\textwidth}
    \centering
    \includegraphics[scale=.25]{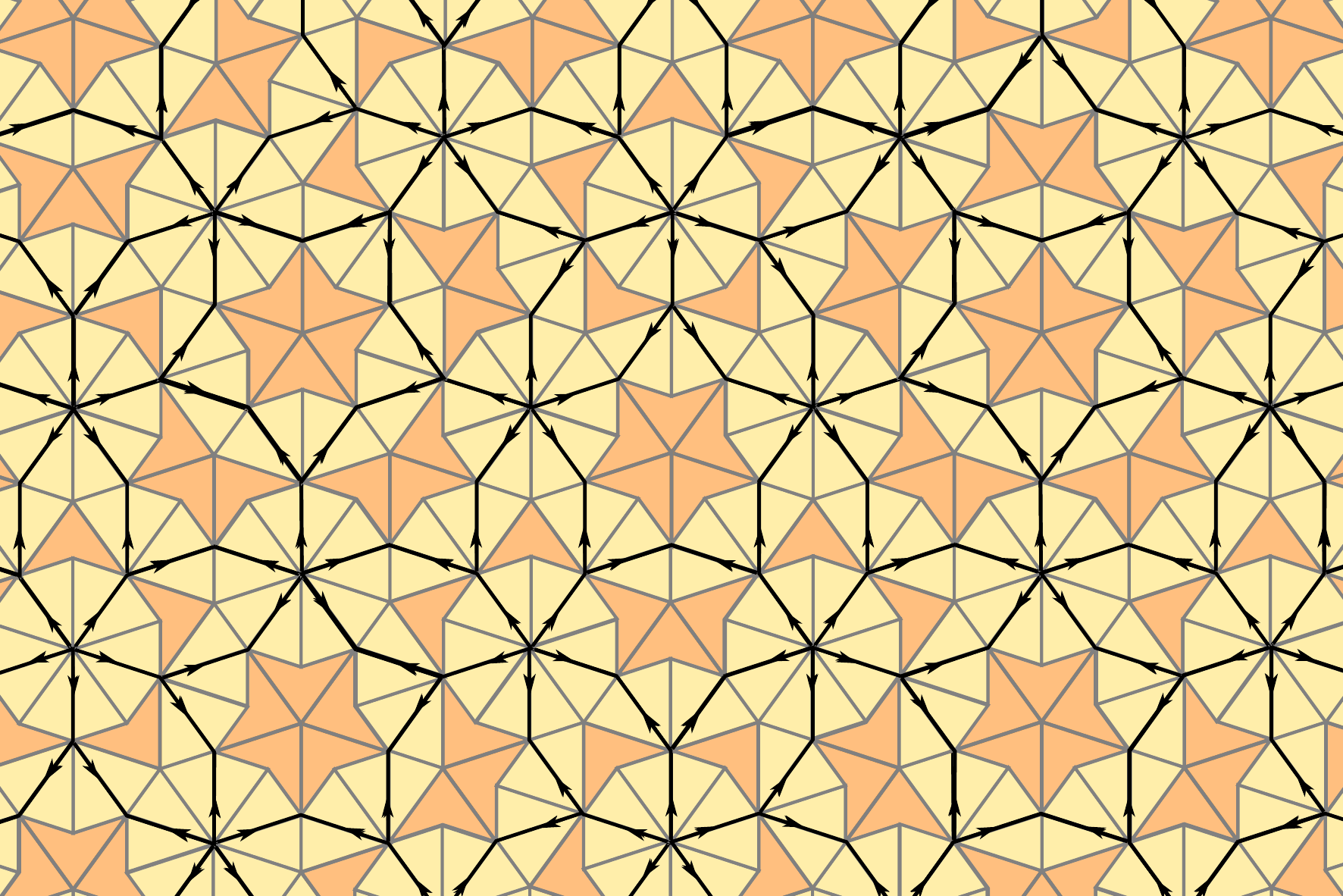}
    \caption{HBS and P2 Penrose tilings superimposed.}\label{fig:hbs-P2}
    \end{subfigure}
    \caption{(a) From a P3 tiling, remove all edges with a double arrow and the vertices they point to. (b) From a P2 tiling, trace the axis of symmetry of each kite, with an arrow pointing to the wide angle, then erase all edges of kites and darts.}\label{fig:hbs-P3-P2}
\end{figure}

HBS tilings are MLD with Penrose tilings of the three types.
We first illustrate the local mappings to and from HBS and Penrose tilings. 
From HBS to Penrose, simply decorate HBS tiles with (parts of) Penrose tiles of the chosen type, always in the same way.
When the Penrose rhombs are marked with arrows as in Figure \ref{P3}, each HBS tile is obtained by composing the rhombs as explained in Figure \ref{fig:hbs-P3}.
HBS tiles are a composition of darts and half-kites as described in Figure \ref{fig:hbs-P2}, with edges of length 1.
With the stars and boats, the resemblance between HBS and P1 tilings is obvious, but the local mapping will be presented further, along with a few additional observations.

\begin{figure}[h]
    \centering
    \begin{subfigure}[b]{0.48\textwidth}
    \centering
    \includegraphics[width=0.6\textwidth]{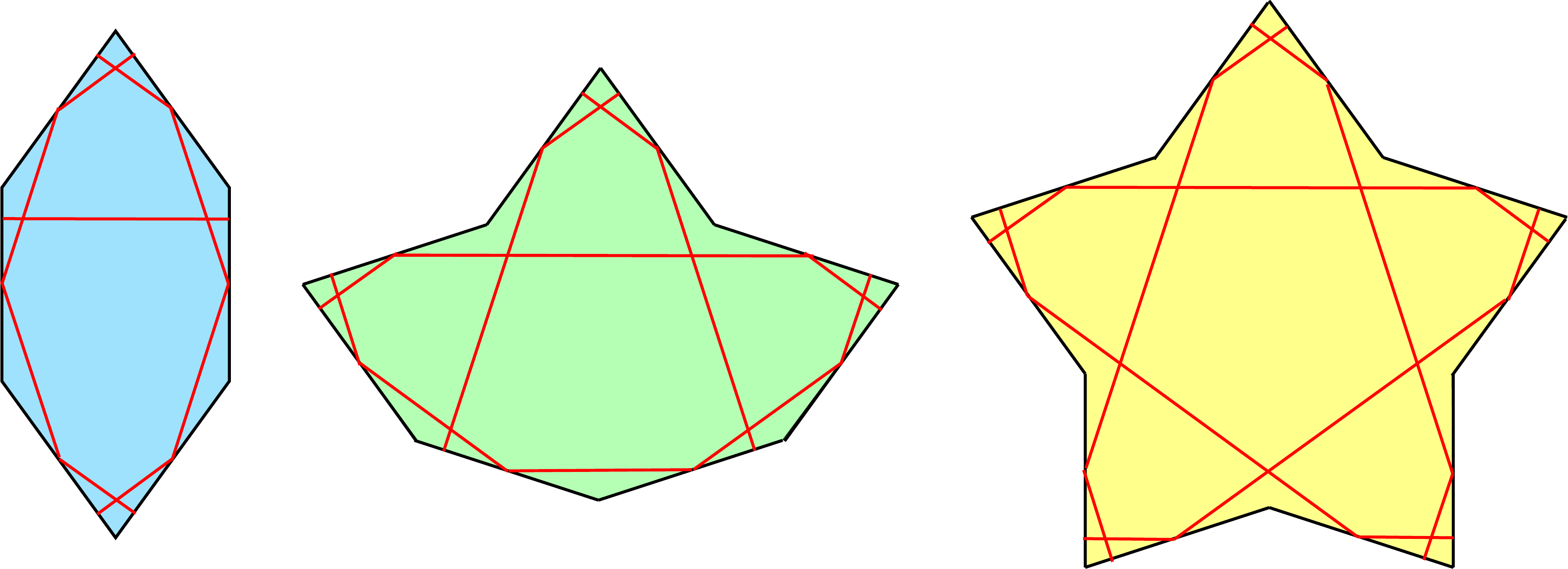}  \caption{HBS tiles with Ammann bar decorations.}
    \label{fig:tileset-BA}
    \end{subfigure}
    \hfill
    \begin{subfigure}[b]{0.48\textwidth}
    \centering
    \includegraphics[width=0.6\textwidth]{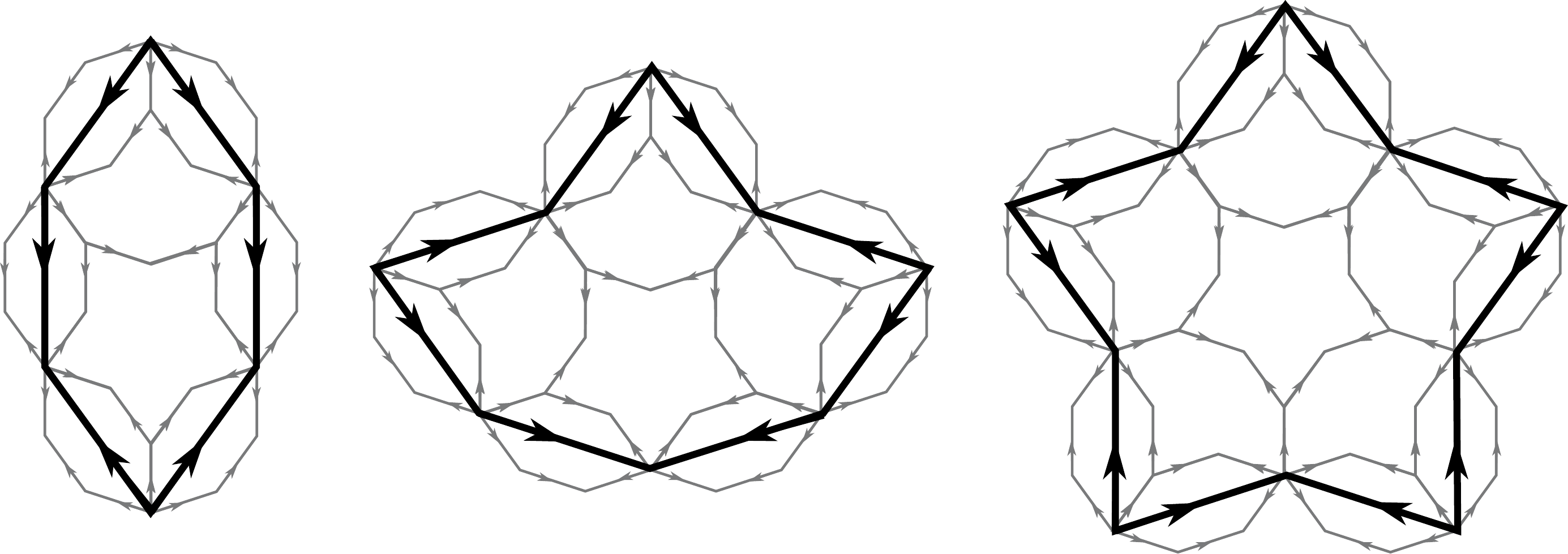}
    \caption{$\varphi^2$-composition of HBS tiles.}
    \label{fig:phi2-composition}
    \end{subfigure}
    \caption{HBS tiles with (a) Ammann segments (arrows are omitted) and (b) their $\varphi^2$-decomposition.}\label{fig:HBS-AB-decomp}
\end{figure}

Henley's paper was cited more than 300 times, but mostly by physicists who experimented different arrangements of atoms based on decorations on HBS shapes -- and not always with the same assembly rules.
Steurer authored a consistent survey of structure research in quasicrystals \cite{Steurer2004} comparing such studies, followed with a book on quasicrystals with Deloudi including a lot more theoretical content \cite{SD2009}. 
In particular, Fig. 1.7 page 24 illustrates the decorations of HBS tiles by Ammann line segments as in Figure \ref{fig:tileset-BA} and give the relative vertex frequencies of P3 tilings, including the conﬁgurations which transform into HBS tiles.
The ratio of hexagons to boats to stars is $\sqrt{5}\varphi : \sqrt{5} : 1$, so that the frequencies of the tiles, as computed by Olamy and Kléman \cite{OK1989}, are 
$$f_H=7-4\varphi=\sqrt{5}\varphi^{-3}\simeq52,8\% ~,\quad 
f_B=7\varphi-11=\sqrt{5}\varphi^{-4}\simeq32,6\% ~,\quad 
f_S=5-3\varphi=\varphi^{-4}\simeq14,6\%.$$
Gummelt \cite{Gummelt2006} also gives the $\varphi^2$-composition of HBS tilings, drawn in Figure \ref{fig:phi2-composition} with the arrows. 
She came up with a decagon covering model, equivalent to the HBS tiling model, as did Lück earlier \cite{LUCK1990}.
Their models were recently compared by Steurer \cite{Steurer2021}.


\section{The Star tileset}

The tileset in Figure \ref{star-tileset-kd} yields the same tilings as the HBS tiles once the decorations are removed.
It was derived from P2 tilings, but not in the same way as in Figure \ref{fig:hbs-P2}.
The construction is still quite simple and in addition to the $\varphi^2$-composition mentioned above, we have a $\varphi$-composition, which is what interested us for the construction in \cite{porrier2023}.
Vertex colors are strongly related to vertex configurations in P2, hence to vertex configurations in HBS tilings when you look at the compositions.

\begin{figure}[ht!]
    \centering
    \includegraphics[width=0.7\textwidth]{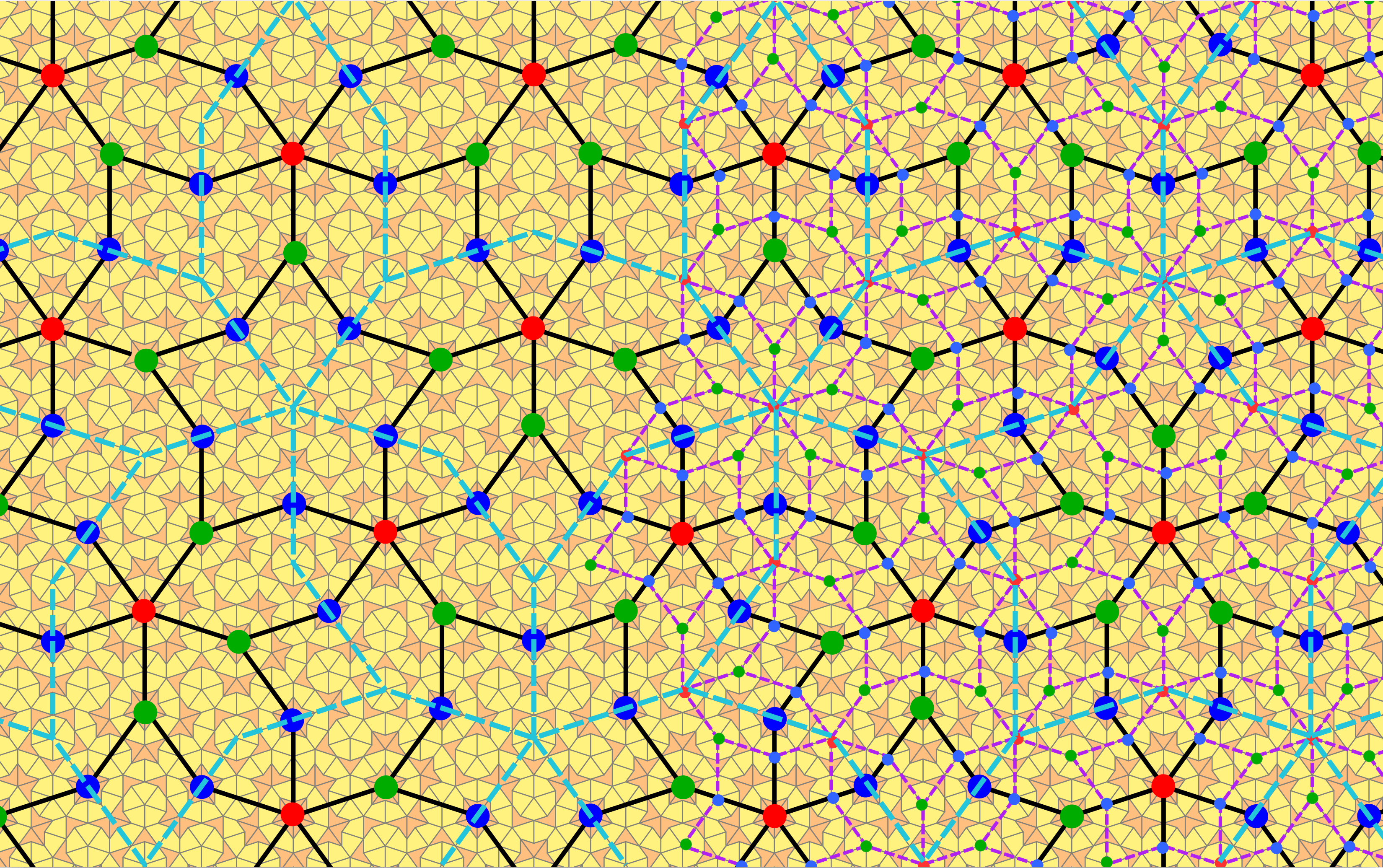}
    \caption{P2 and Star tilings superimposed.}
    \label{fig:star-tiling-example}
\end{figure}

As can be seen in Figure \ref{fig:hbs-P2} and considering the local rules of P2 tilings, the shortest distance between the centers of two stars\footnote{``Star'' vertex configurations in a P2 tiling or stars in a HBS tiling.} is $3+2\varphi^{-1}=\varphi^3$, which occurs when two HBS-stars are incident to the same HBS-edge (or kite of the underlying PT).
If you join the centers of the stars, the HBS shapes appear again but larger, composed with many kites (or half-kites) and darts as in Figure \ref{fig:star-tiling-example}.
The vertices are colored according to their degree, which corresponds to the type of P2-star they lie on: the star can be tangent to either 0 (red), 1 (green) or 2 (blue) suns.
If you put on each edge an arrow pointing to the greatest number (among the endpoints), you get exactly the HBS tiles as in Figure \ref{fig:hbs-tiles}.
Yet here there is more information encoded in the vertices, hence more local constraints.
The hexagons and boats all have exactly the same decorations but the stars come in three versions in what we call the Star tileset.\\

\begin{figure}[h]
    \centering
    \includegraphics[width=0.6\textwidth]{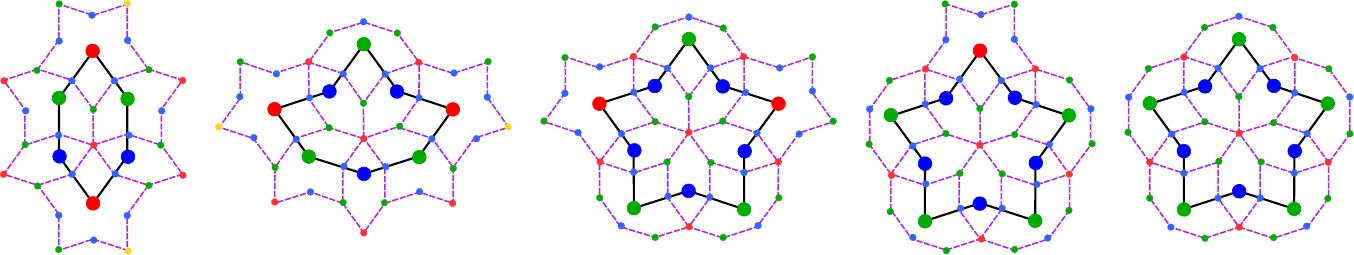}
    \caption{Star tileset decomposition.}
    \label{fig:star-decomposition}
\end{figure}

Obviously, the Ammann bar decorations of HBS tiles are still valid for the same shapes in the Star tileset.
Each vertex of the star tiling lies inside a small polygon formed by $n$ Ammann bars.
The label of the vertex is then $5-n$ (or the corresponding color).
Since the vertices of the HBS tiling are the stars of the P2 tiling, and a star is obtained by inflating a sun, joining the suns which are at distance $2+\varphi^{-1}=\varphi^2$ (the minimum) as in Figure \ref{fig:star-tiling-example} (magenta dashed lines) yields the same HBS tiling as first inflating once the P2 tiling and then joining the stars.
The vertex colors are then given by the number of queen configurations intersecting with the sun.
The resulting $\varphi$-decomposition of each tile is resumed in Figure \ref{fig:star-decomposition}.
Apply it twice to get the $\varphi^2$-decomposition mentioned above.
Now if you apply the $\varphi$-decomposition a third time and decorate the smaller shapes as in Figure \ref{fig:hbs-P2} (usual decoration of HBS tiles with half-kites and darts), the initial Star tiles will have the same decorations as in Figure \ref{star-tileset-kd}
.
\begin{figure}[h!]
    \centering
    \begin{subfigure}[b]{0.48\textwidth}
    \centering
    \includegraphics[scale=.25]{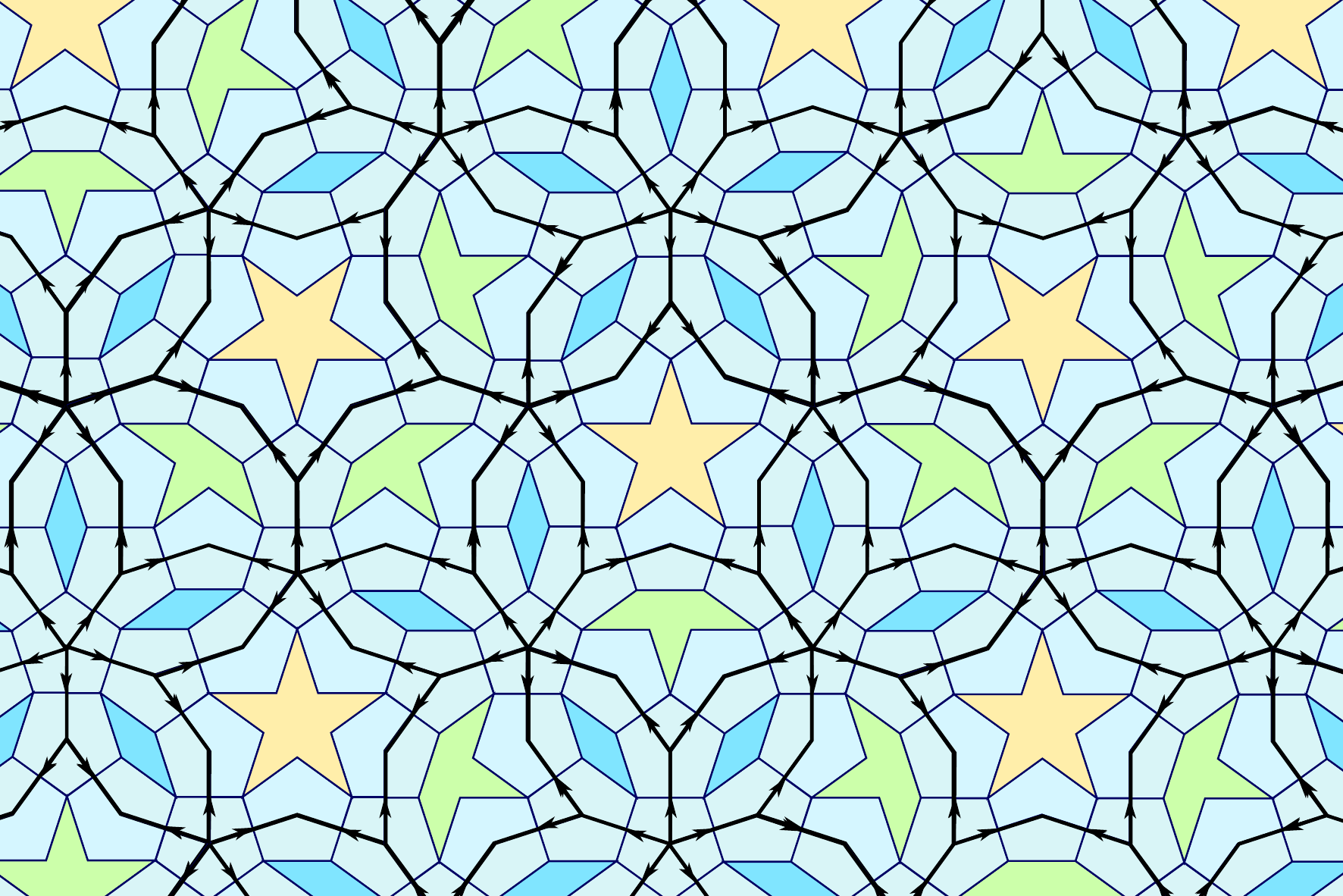}
    \caption{HBS and P1 Penrose tilings superimposed. }\label{fig:hbs-P1}
    \end{subfigure}
    \hfill
    \begin{subfigure}[b]{0.48\textwidth}
    \centering
    \includegraphics[scale=.25]{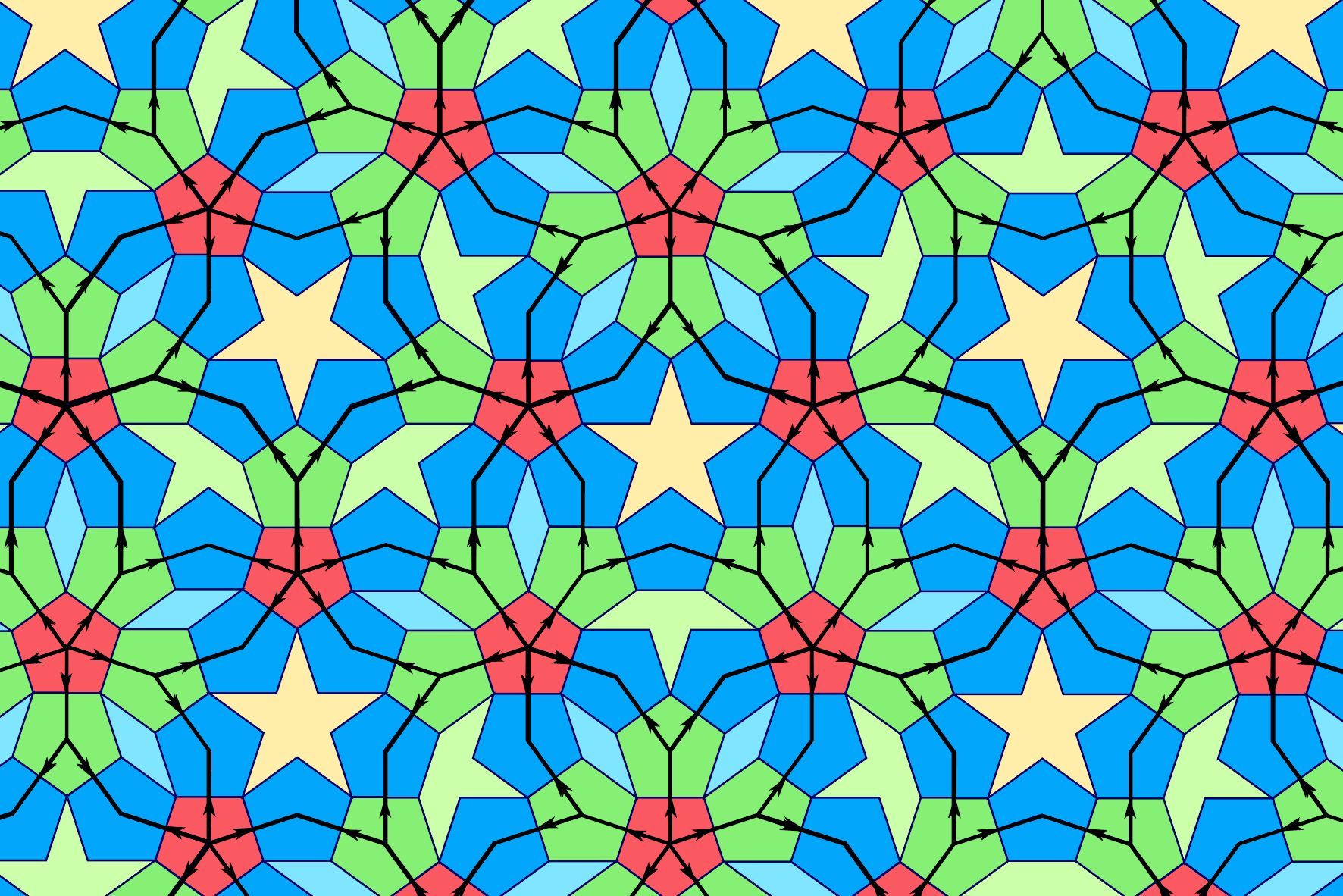}
    \caption[Star and P1 Penrose tilings superimposed.]{Star and P1 Penrose tilings superimposed.}\label{fig:star-P1-example}
    \end{subfigure}
    \caption{Each P1-star is inside a HBS-star, each P1-boat inside a HBS-boat, and each P1-diamond inside a HBS-hexagon, while all pentagons are decorated with the edges of HBS tiles.
    In the Star tiling, pentagons are colored according to their type (Figure \ref{fig:P1-tiles2}).}
\end{figure}
\begin{figure}[h!]
    \centering
    \begin{subfigure}[b]{0.65\textwidth}
    \centering
    \includegraphics[scale=0.8]{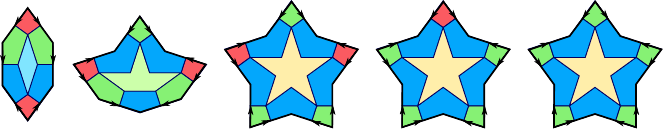}
    \caption{Star tileset with P1 decorations.}\label{star-P1-tiles}
    \end{subfigure}
    \hfill
    \begin{subfigure}[b]{0.3\textwidth}
    \centering
    \includegraphics[scale=0.8]{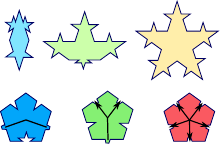}
    \caption{P1 with colored pentagons.}
    \label{fig:P1-tiles2}
    \end{subfigure}
    \caption{Duality between P1 and Star tilesets.}\label{fig:P1-star}
\end{figure}
Equivalently, from an HBS tiling with the tiles decorated as in Figure \ref{fig:hbs-P2} you can keep the HBS shapes as is, apply three $\varphi$-decompositions to kites and darts, and then appropriately place the vertex colors on hexagons and boats.
Actually, you can even place the colors from the start, considering how vertex configurations in P2 tilings substitute.
As for the substitution with P1, the label of a pentagon is the number of HBS-arrows pointing to it.
Thus each type of pentagon from the original P1 tileset can simply get the corresponding color, as illustrated in Figures \ref{fig:star-P1-example} and \ref{fig:P1-star}, and if we add the arrows on P1 tiles then the corresponding HBS tiling appears on any P1 tiling.

\begin{figure}[h]
    \centering
    \includegraphics[width=\textwidth]{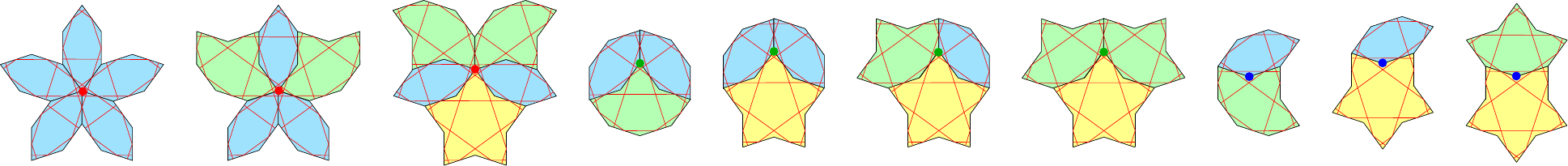}
    \caption[Vertex neighborhoods in HBS and Star tilings.]{Vertex neighborhoods in HBS and Star tilings, with different possible colors on some vertices of the stars. When stars, boats and hexagons are decomposed, we obtain one of those centered on a red vertex (from left to right), which we call respectively bellflower, orchid and pansy.}\label{fig:star-config}
\end{figure}

HBS vertex configurations are given in Figure \ref{fig:star-config}.
The main difference between HBS and Star tilings is the amount of information which can be deduced from a local configuration. 
As you can see in Figure \ref{fig:star-decomposition}, 
each vertex color in a Star tiling decomposes into an hexagon, a boat or a star according to its color.
Hence the ratio of blue to green to red is the same as hexagons to boats to stars, that is $\sqrt{5}\varphi : \sqrt{5} : 1$.
Since the frequency of stars in HBS tilings is $f_S=\varphi^{-4}$, we have
$$f_H=\sqrt{5}\varphi^{-3} ~,\quad 
f_B=\sqrt{5}\varphi^{-4} ~,\quad 
f_2=\sqrt{5}\varphi^{-7}\simeq7,7\% ~,\quad 
f_1=\sqrt{5}\varphi^{-8}\simeq4,76\% ~,\quad 
f_0=\varphi^{-8}\simeq2,13\%$$
where $f_i$ is the frequency of the star S$_i$ with $i$ red vertices.
These frequencies can also be computed using those of vertex configurations in Penrose tilings.

\begin{figure}[h]
    \centering
    \begin{subfigure}[b]{0.3\textwidth}
    \centering
    \includegraphics[scale=0.5]{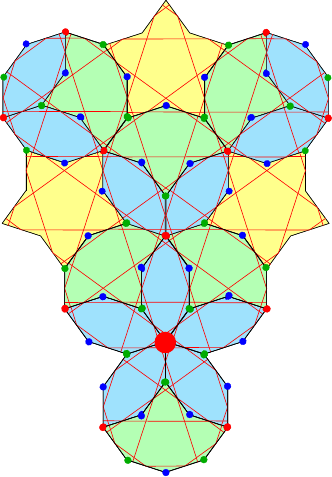}
    \caption{Kingdom of the bellflower.}\label{fig:bellflower}
    \end{subfigure}
    \hfill
    \begin{subfigure}[b]{0.4\textwidth}
    \centering
    \includegraphics[scale=0.5]{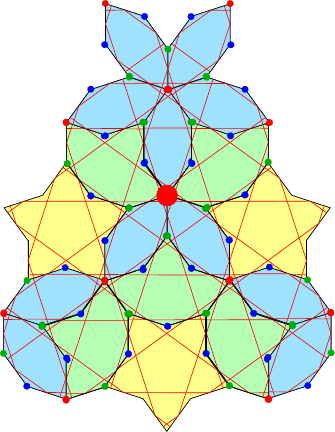}
    \caption{Kingdom of the orchid.}
    \label{fig:orchid}
    \end{subfigure}
\end{figure}
\begin{figure*}[h]
    \centering
    \vspace{-1.8cm}
    \begin{subfigure}[b]{0.2\textwidth}
    \centering
    \includegraphics[scale=0.5]{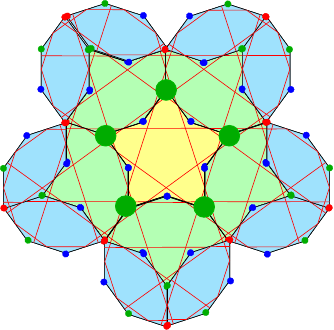}
    \caption*{(c) Kingdom of S0.}\label{fig:star0}
    \end{subfigure}
    \hfill
    \begin{subfigure}[b]{0.44\textwidth}
    \centering
    \includegraphics[scale=0.5]{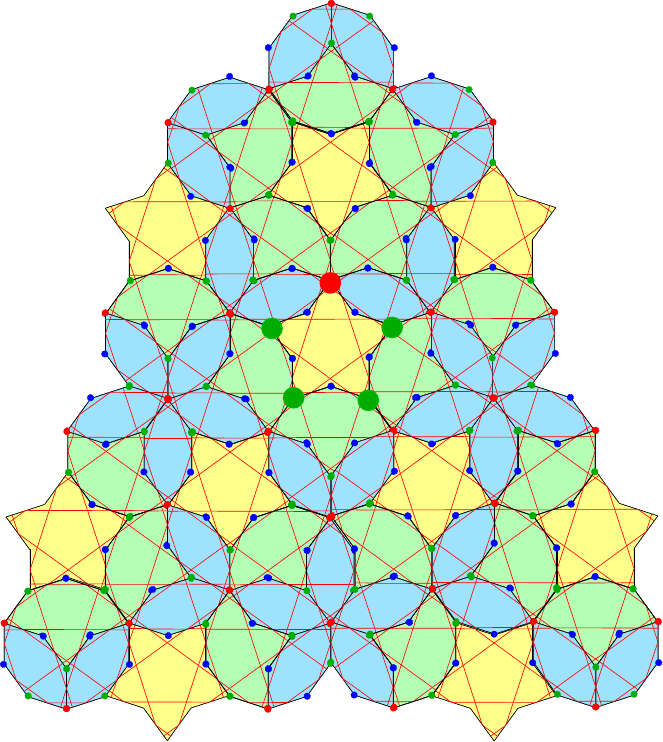}
    \caption*{(d) Kingdom of S1.}
    \label{fig:star1}
    \end{subfigure}
    \hfill
    \begin{subfigure}[b]{0.32\textwidth}
    \centering
    \includegraphics[scale=0.5]{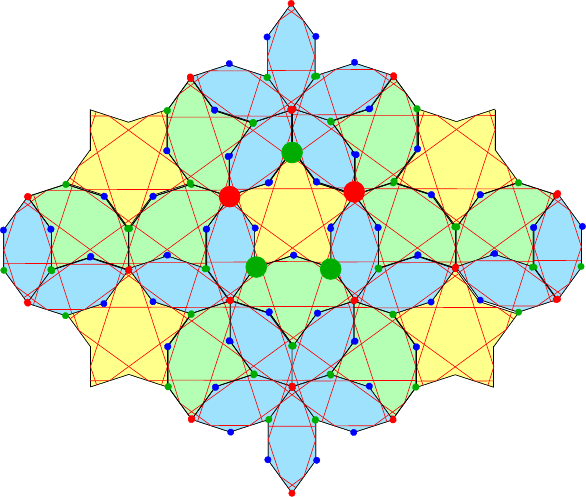}
    \caption*{(e) Kingdom of S2.}
    \label{fig:star2}
    \end{subfigure}
    \caption{Kingdoms of the bellflower, the orchid and the 3 stars.}\label{fig:stars-kingdoms}
\end{figure*}

As in Penrose tilings, each vertex configuration can force a whole set of tiles (in the tiling) which is called \textit{empire}.
This means that anytime a given vertex configuration appears in any tiling of the same type, 
all the tiles forming the empire are placed exactly in the same way relatively to it.
Since several empires are disconnected and can even contain infinitely many tiles, we only represent here the largest connected component, which includes the vertex configuration. We call it \textit{kingdom} (or \textit{local empire}).
For the bellflower and the orchid, since they contain no stars, vertex colors could be omitted (as in HBS tilings), but when a vertex configuration includes a star, the colors have a significant impact, as you can see in Figure \ref{fig:stars-kingdoms}.
In particular, the star S1 forces a star S0 above and two stars S2 below.
These three kingdoms resemble those of respectively the star, the king and the queen in P2 tilings.


\section{Gemstones tileset}

\begin{figure}[h!]
    \centering
    \includegraphics[scale=0.25]{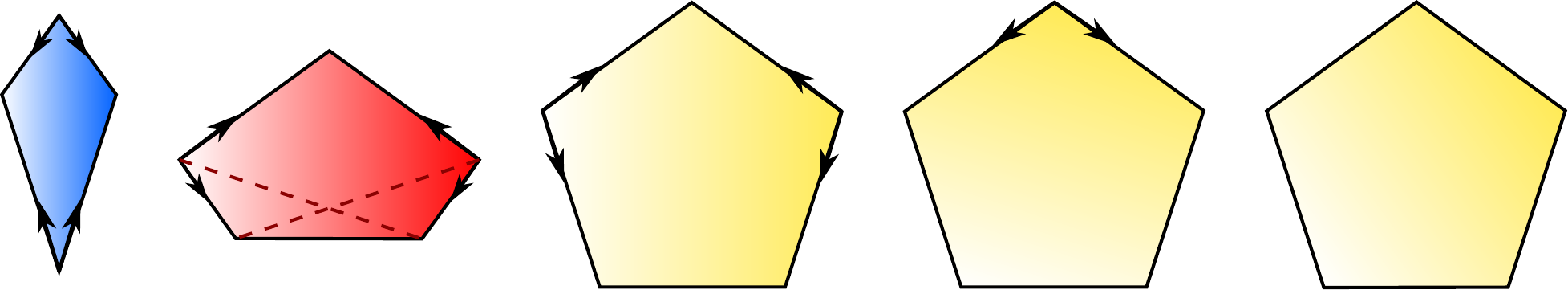}
    \caption{Gemstones tileset: sapphire, ruby and topazes. 
    Edges have length 1 and $2\sin{\frac{2\pi}{5}}$.
    The intersection point of the dashed lines is the ``center'' of the ruby for the $\varphi$-composition in Figure \ref{fig:gemstones-subst}.}
    \label{fig:gemstones-tileset}
\end{figure}

The Star tileset can be modified to get convex tiles with only two labels at their vertices, observing that a blue vertex of an hexagon, a boat or a star is actually not a vertex of the Star tiling.
Thus we obtain the \textit{Gemstones} tileset in Figure \ref{fig:gemstones-tileset}.
On Figure \ref{fig:gemstones-example}, one can see how HBS tiles are deformed to get gemstones.
Substitution from P3 is also quite simple.
Trace the long diagonal of each thin rhomb.
For each Q configuration (hexagon), compose the fat rhomb with its two adjacent half thin rhombs.
Everywhere else just erase all edges of the P3 tiling.
Conversely, Gemstones are easily decorated with rhombuses (including arrows or other matching rules).
As for P2, the suns and jacks are vertices of Gemstones: trace a (short) segment between those vertices whenever they share a kite, and a (long) segment when they are separated by two kites which share a short edge.
\begin{figure}[h]
    \centering
    \includegraphics[width=0.7\textwidth]{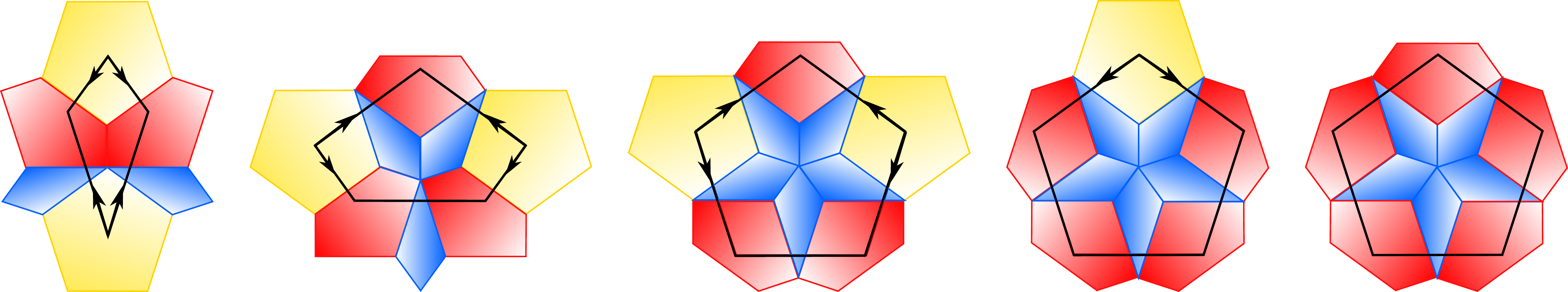}
    \caption{$\varphi$-decomposition of Gemstones.}\label{fig:gemstones-subst}
\end{figure}
\begin{figure}[h!]
    \centering
    \includegraphics[scale=.25]{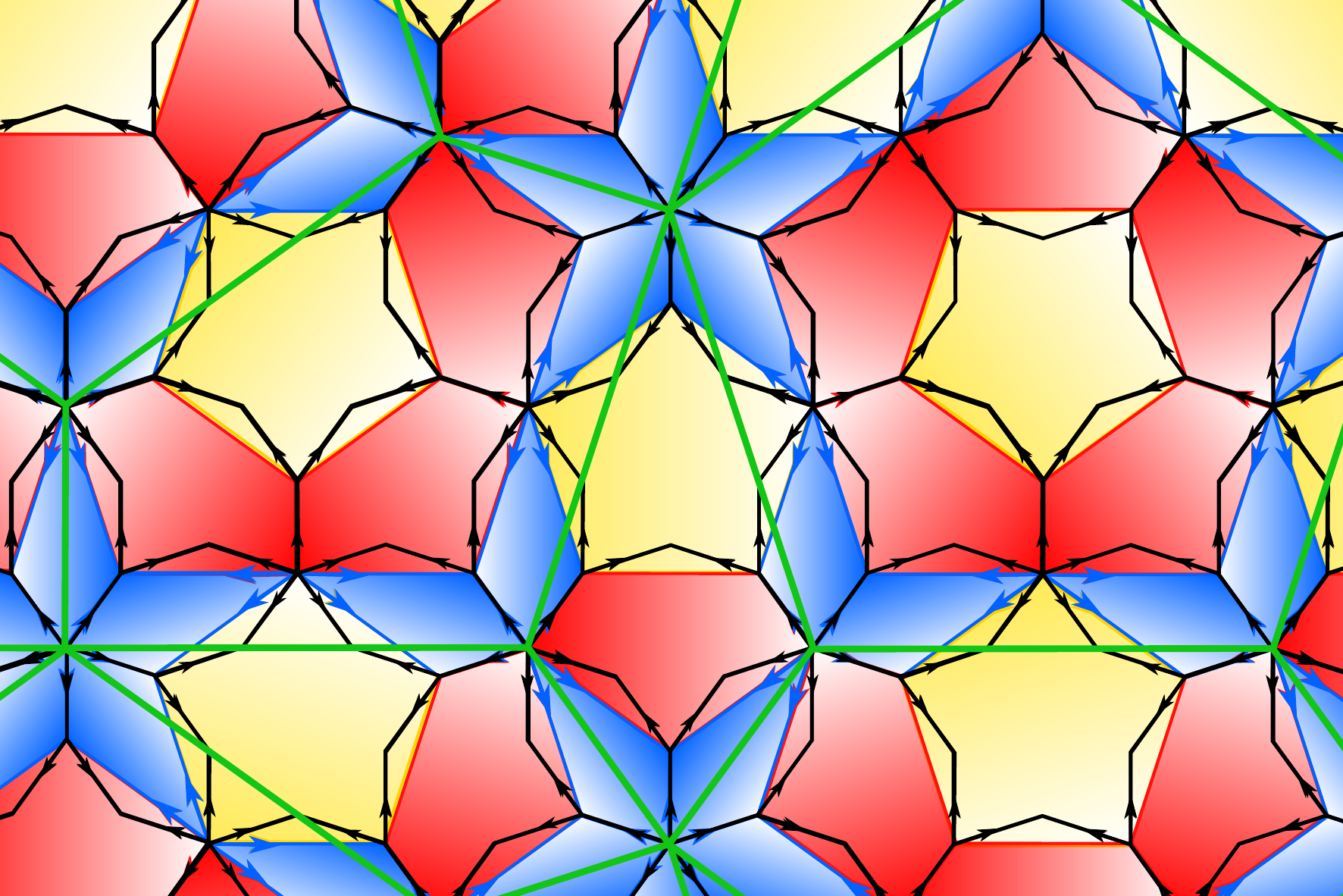}
    \caption{Gemstones and HBS tilings superimposed, along with the $\varphi^2$-composition in green.}\label{fig:gemstones-example}
\end{figure}


\bibliographystyle{eptcs}
\bibliography{hbs}

\begin{thebibliography}{10}
\providecommand{\bibitemdeclare}[2]{}
\providecommand{\surnamestart}{}
\providecommand{\surnameend}{}
\providecommand{\urlprefix}{Available at }
\providecommand{\url}[1]{\texttt{#1}}
\providecommand{\href}[2]{\texttt{#2}}
\providecommand{\urlalt}[2]{\href{#1}{#2}}
\providecommand{\doi}[1]{doi:\urlalt{https://doi.org/#1}{#1}}
\providecommand{\eprint}[1]{arXiv:\urlalt{https://arxiv.org/abs/#1}{#1}}
\providecommand{\bibinfo}[2]{#2}

\bibitemdeclare{book}{BaakeGrimm2013}
\bibitem{BaakeGrimm2013}
\bibinfo{author}{Michael \surnamestart Baake\surnameend} \&
  \bibinfo{author}{Uwe \surnamestart Grimm\surnameend} (\bibinfo{year}{2013}):
  \emph{\bibinfo{title}{Aperiodic Order. Vol 1. A Mathematical Invitation}}.
\newblock {\slshape \bibinfo{series}{Encyclopedia of Mathematics and its
  Applications}} \bibinfo{volume}{149}, \bibinfo{publisher}{Cambridge
  University Press}, \bibinfo{address}{Cambridge}.
\newblock \urlprefix\url{http://oro.open.ac.uk/38933/}.

\bibitemdeclare{article}{gardner1977}
\bibitem{gardner1977}
\bibinfo{author}{Martin \surnamestart Gardner\surnameend}
  (\bibinfo{year}{1977}): \emph{\bibinfo{title}{MATHEMATICAL GAMES}}.
\newblock {\slshape \bibinfo{journal}{Scientific American}}
  \bibinfo{volume}{236}(\bibinfo{number}{1}), pp. \bibinfo{pages}{110--121},
  \doi{10.1038/scientificamerican0177-110}.

\bibitemdeclare{book}{grunbaum1987}
\bibitem{grunbaum1987}
\bibinfo{author}{Branko \surnamestart Grünbaum\surnameend} \&
  \bibinfo{author}{G.~C. \surnamestart Shephard\surnameend}
  (\bibinfo{year}{1987}): \emph{\bibinfo{title}{Tilings and Patterns}}.
\newblock \bibinfo{publisher}{W. H. Freeman and Company, New York}.
\newblock \bibinfo{note}{ISBN 0-7167-1193-1}.

\bibitemdeclare{article}{Gummelt2006}
\bibitem{Gummelt2006}
\bibinfo{author}{Petra \surnamestart Gummelt\surnameend}
  (\bibinfo{year}{2006}): \emph{\bibinfo{title}{Decagon covering model and
  equivalent {HBS}-tiling model}}.
\newblock {\slshape \bibinfo{journal}{Zeitschrift für Kristallographie -
  Crystalline Materials}} \bibinfo{volume}{221}(\bibinfo{number}{8}), pp.
  \bibinfo{pages}{582--588}, \doi{10.1524/zkri.2006.221.8.582}.

\bibitemdeclare{article}{Henley1986}
\bibitem{Henley1986}
\bibinfo{author}{Christopher~L. \surnamestart Henley\surnameend}
  (\bibinfo{year}{1986}): \emph{\bibinfo{title}{{Sphere packings and local
  environments in Penrose tilings}}}.
\newblock {\slshape \bibinfo{journal}{Phys. Rev. B}} \bibinfo{volume}{34}, pp.
  \bibinfo{pages}{797--816}, \doi{10.1103/PhysRevB.34.797}.
\newblock \bibinfo{note}{Reprinted in \cite{steinhardt1987}}.

\bibitemdeclare{article}{LUCK1990}
\bibitem{LUCK1990}
\bibinfo{author}{Reinhard \surnamestart Lück\surnameend}
  (\bibinfo{year}{1990}): \emph{\bibinfo{title}{Penrose sublattices}}.
\newblock {\slshape \bibinfo{journal}{Journal of Non-Crystalline Solids}}
  \bibinfo{volume}{117-118}, pp. \bibinfo{pages}{832--835},
  \doi{10.1016/0022-3093(90)90657-8}.

\bibitemdeclare{article}{OK1989}
\bibitem{OK1989}
\bibinfo{author}{\surnamestart {Olamy, Z.}\surnameend} \&
  \bibinfo{author}{\surnamestart {Kl\'eman, M.}\surnameend}
  (\bibinfo{year}{1989}): \emph{\bibinfo{title}{A two dimensional aperiodic
  dense tiling}}.
\newblock {\slshape \bibinfo{journal}{J. Phys. France}}
  \bibinfo{volume}{50}(\bibinfo{number}{1}), pp. \bibinfo{pages}{19--33},
  \doi{10.1051/jphys:0198900500101900}.

\bibitemdeclare{article}{penrose1978}
\bibitem{penrose1978}
\bibinfo{author}{Roger \surnamestart Penrose\surnameend}
  (\bibinfo{year}{1978}): \emph{\bibinfo{title}{{Pentaplexity}}}.
\newblock {\slshape \bibinfo{journal}{Math. Intelligencer}}
  \bibinfo{volume}{2}(\bibinfo{number}{1}), pp. \bibinfo{pages}{32--37},
  \doi{10.1007/BF03024384}.

\bibitemdeclare{article}{porrier2020}
\bibitem{porrier2020}
\bibinfo{author}{Carole \surnamestart Porrier\surnameend} \&
  \bibinfo{author}{Alexandre \surnamestart Blondin~Massé\surnameend}
  (\bibinfo{year}{2020}): \emph{\bibinfo{title}{The {L}eaf Function of Graphs
  Associated with {P}enrose Tilings}}.
\newblock {\slshape \bibinfo{journal}{International Journal of Graph
  Computing}} \bibinfo{volume}{1}, pp. \bibinfo{pages}{1--24},
  \doi{10.35708/GC1868-126721}.

\bibitemdeclare{misc}{porrier2023}
\bibitem{porrier2023}
\bibinfo{author}{Carole \surnamestart Porrier\surnameend},
  \bibinfo{author}{Alain \surnamestart Goupil\surnameend} \&
  \bibinfo{author}{Alexandre \surnamestart Blondin~Massé\surnameend}
  (\bibinfo{year}{2023}): \emph{\bibinfo{title}{The Leaf Function of Penrose P2
  Graphs}}.
\newblock \eprint{2312.08262}.

\bibitemdeclare{book}{Senechal1995}
\bibitem{Senechal1995}
\bibinfo{author}{M.~\surnamestart Senechal\surnameend} (\bibinfo{year}{1995}):
  \emph{\bibinfo{title}{{Quasicrystals and Geometry}}}.
\newblock \bibinfo{publisher}{Cambridge University Press}.
\newblock \bibinfo{note}{ISBN 0521372593}.

\bibitemdeclare{book}{steinhardt1987}
\bibitem{steinhardt1987}
\bibinfo{author}{Paul~J \surnamestart Steinhardt\surnameend} \&
  \bibinfo{author}{Stellan \surnamestart Ostlund\surnameend}
  (\bibinfo{year}{1987}): \emph{\bibinfo{title}{The Physics of Quasicrystals}}.
\newblock \bibinfo{publisher}{World Scientific}, \doi{10.1142/0391}.
\newblock \bibinfo{note}{Collection of reprints.}

\bibitemdeclare{article}{Steurer2004}
\bibitem{Steurer2004}
\bibinfo{author}{Walter \surnamestart Steurer\surnameend}
  (\bibinfo{year}{2004}): \emph{\bibinfo{title}{{Twenty years of structure
  research on quasicrystals. Part I. Pentagonal, octagonal, decagonal and
  dodecagonal quasicrystals}}}.
\newblock {\slshape \bibinfo{journal}{Zeitschrift für Kristallographie -
  Crystalline Materials}} \bibinfo{volume}{219}(\bibinfo{number}{7}), pp.
  \bibinfo{pages}{391--446}, \doi{10.1524/zkri.219.7.391.35643}.

\bibitemdeclare{article}{Steurer2021}
\bibitem{Steurer2021}
\bibinfo{author}{Walter \surnamestart Steurer\surnameend}
  (\bibinfo{year}{2021}): \emph{\bibinfo{title}{{Gummelt versus Lück decagon
  covering and beyond. Implications for decagonal quasicrystals}}}.
\newblock {\slshape \bibinfo{journal}{Acta Crystallographica Section A}}
  \bibinfo{volume}{77}(\bibinfo{number}{1}), pp. \bibinfo{pages}{36--41},
  \doi{10.1107/S2053273320015181}.

\bibitemdeclare{book}{SD2009}
\bibitem{SD2009}
\bibinfo{author}{Walter \surnamestart Steurer\surnameend} \&
  \bibinfo{author}{Sofia \surnamestart Deloudi\surnameend}
  (\bibinfo{year}{2009}): \emph{\bibinfo{title}{Crystallography of
  Quasicrystals: Concepts, Methods and Structures}}, \bibinfo{edition}{1}
  edition.
\newblock \bibinfo{series}{Springer Series in Materials Science №126},
  \bibinfo{publisher}{Springer}, \doi{10.1007/978-3-642-01899-2}.

\end{thebibliography}

\end{document}